\magnification 1200
\input amstex
\documentstyle{amsppt}
\loadbold
\NoBlackBoxes

\def\ju{\vskip .8truecm plus .1truecm minus .1truecm}

\def\sju{\vskip .4truecm plus .1truecm minus .1truecm}
\def\ets{\emptyset}
\def\se{\subseteq}
\def\al{\alpha}

\def\ka{\kappa}
\def\be{\beta}

\def\om{\omega}
\def\sst#1#2{\{ {#1}:{#2} \}}
\def\lra{\longrightarrow}
\def\rmod{R\bold{Mod}}
\def\shimply{\Rightarrow} 
\def\pr{^\prime}
\def\ann{\text{Ann}}
\def\lla{\longleftarrow}
\def\hookr{\hookrightarrow}
\def\op{\oplus}
\def\spe{\supseteq}
\def\homc#1#2#3{\hbox{\rm Hom}_{#1}({#2},{#3})}
\def\ima{\hbox{\rm Im\,}}
\def\ese#1{\buildrel{#1}\over\lra}
\def\ov{\overline}

\document

\font\titlefontd=cmssdc10 at 20pt

\font\ninesl=ptmri7t at 15pt 
\font\tencmmib=cmmib10 \skewchar\tencmmib='177   

\nologo
\def\noi{\noindent}
\def\prf{\noi{\it Proof.\quad}}
\def\qed{\hfill$\square$}

\hsize 15.5truecm
\vsize 20truecm
\def\narrower{\advance\leftskip by 2.5truecm 
                \advance\rightskip by 2.5truecm}

\def\csju{\vskip 0.3truecm}
\def\e{\bold e}

\vskip 3truecm

\centerline{\titlefontd Kappa-Slender Modules\footnote[^*]{The first (handwritten, unpublished) version of this paper dates back to 1985.}}
\rightheadtext\nofrills{\it $\ka$-Slender Modules}
\leftheadtext\nofrills{\it $\ka$-Slender Modules}

\vskip 0.5truecm
\centerline{By}
\vskip 0.5truecm
\centerline{\ninesl Radoslav Dimitric\footnote[]{email: dimitricr at member.ams.org}}
\vskip 1truecm

\noi{\bf Abstract.} For an arbitrary infinite cardinal $\ka$, we define classes of  $c\ka$-slender and $t\ka$-slender modules as well as related classes of $h\ka$-modules and initiate a study of these classes. 
\ju

\noi{\bf 0. Preliminaries} 
\csju
The Axiom of Choice (AC) is assumed which in particular implies that the class of ordinals may be well-ordered and that the cardinals may be then enumerated by ordinals and replaced by alephs $\aleph_\alpha$; the smallest ordinal of cardinality $\aleph_\al$ is denoted by $\om_\al$ and $\om=\om_0$ will be used to denote the first infinite ordinal. We will freely interchange $\aleph_\al$ and $\om_\al$. Given ordinals $\be$ and $\al$, we say that $\al$ is {\it the cofinality} of $\be$ (notation $\al= cf\be$), if $\al$ is the smallest ordinal order-isomorphic to a subset $A$ of $\be$ such that the least upper bound of $A$ is $\be$. An ordinal $\be$ is said to be {\it regular}, if $cf\,\be=\be$; otherwise, $\be$ is a {\it singular ordinal}. $\al=\be+1$ is called the {\it successor ordinal} to $\be$. If $\al$ cannot be represented in this form, then $\al$ is said to be a {\it limit ordinal}. A cardinal $\lambda$ is called {\it regular} if it is not the supremum of $<\lambda$ cardinals which are $<\lambda$; equivalently $\lambda$ is not the sum of fewer than $\lambda$ cardinals which are $<\lambda$, equivalently, if $cf\lambda=\lambda$; other cardinals are called {\it singular cardinals}. $\aleph_n$'s are regular cardinals, while $\aleph_\om$ is singular, since $\sum_{n<\om}\aleph_n=\aleph_\om$. The first cardinal greater than cardinal $\ka$ is denoted by $\ka^+$ and cardinals of this form are called {\it successor cardinals}; other cardinals are called {\it limit cardinals}. Every successor cardinal is regular.

For a non-empty infinite set $I$, $\Cal I\se 2^I$ is an {\it ideal} on $I$, if $A\in\Cal I$ and $B\se A$ imply that $B\in\Cal I$ and any finite union of members of $\Cal I$ is again in $\Cal I$. This is a proper ideal, if $I\notin\Cal I$; otherwise, it is improper. Dually, we arrive at the notion of a filter: $\Cal F\se 2^I$ is a {\it filter} on $I$, if $A\in\Cal F$ and $A\se B$ imply that $B\in\Cal F$ and any finite intersection of members of $\Cal F$ is again in $\Cal F$. This is a proper filter, if $\ets\notin\Cal F$, otherwise it is an improper filter. $\sst{A\se I}{|A|<|I|}$ is a proper ideal on $I$ and,  dually, $\sst{A\se I}{|cA|<|I|}$ is a proper filter on $I$. 
For an infinite cardinal $\ka$, 
If $|I|=\ka$, then $\Cal F_{F\ka}=\sst{X\se I}{|cX|<\ka}$ is a proper filter, 
called the $\ka$-{\it Fr\'echet filter}. The filter is still proper if $\ka < |I|$.
If we do not have restrictions on $\ka$ in relation to $|I|$, the filter
$\Cal F_\ka=\sst{X\se I}{|cX|<\ka}$ is a {\it co}-$\kappa$-{\it filter} that may
be improper ($\Cal F=2^I$). 
Given an infinite cardinal $\ka$ a filter $\Cal F$ is called $\ka$-{\it complete},
if, for every non-empty $\Cal F_1\se\Cal F$ with $|\Cal F_1|<\ka$ we have $\cap\Cal F_1\in\Cal F$. 
With this definition, all filters are $\aleph_0$-complete; an $\aleph_1$-complete
filter will be called {\it countably complete}. One can show that $\Cal F$
is principal iff it is $\ka$-complete, for every $\kappa$. A filter is
$\ka$-{\it incomplete}, if it is not $\ka$-complete; in particular, we define {\it countably incomplete} filters.
For $x=(x_i)_{i\in I}\in\prod_{i\in I} A_i$ denote the zero set zero$(x)=\sst{i\in I}{x_i=0}$ and the complement non-zero set supp$(x)=\sst{i\in I}{x_i\neq 0}$. Given a filter $\Cal F$ on a non-empty index set $I$, the $\Cal F$ subproduct is $\Pi(\Cal F)=\prod^\Cal F_{i\in I}M_i=\sst{x\in\prod M_i}{\text{ zero}(x)\in\Cal F}$. If $\Cal F$ is a co-$\ka$-filter, then the corresponding subproduct will be denoted by $\prod_{i\in I}^\ka M_i$; namely it consists of all the vectors $x=(x_i)_{i\in I}$ with support of cardinality $<\ka$. 

The coordinate vectors $\e_i:I\lra R$ are defined by $\e_i(j)=0$, if $i\neq j$ and $\e_i(i)=1\in R$.

We will work within a category of (left) $R$-modules $\rmod$ and our functions will be morphisms in that category; in particular, $R$ is seen as an object in this category, not in the category of rings. To simplify the discussion, we may assume, if need be,  that rings $R$ are domains with unities and that all modules are unitary. We may assume that cardinalities of our index sets are less than some large cardinal $\frak m$, such as a measurable cardinal, or an inaccessible cardinal. 
\ju

\noi{\bf 1. Coordinate slenderness}
\csju
\noi{\bf Definition 1.} Given an arbitrary (infinite) cardinal $\ka$, define a (left) $R$-module $M$ to be $\ka$-{\it slender}, if, for every index set $I$ of cardinality $\ka$, every family of $R$-modules $A_i$, $i\in I$, and every morphism $$f:\prod_{i\in I}A_i\lra M,\quad  |\sst{i\in I}{f|A_i\neq 0}|<\ka.$$
\csju

\noi In this way, the well-known notion of slender module is a special case, namely of an $\omega_0$-slender module (see Dimitric (2017), for a thorough study of classes of slender objects). Another, more appropriate name we will use is {\it coordinatewise} or {\it c}$\ka$ {\it slenderness}. 

The purpose of this note is to look into $\ka$-slenderness for uncountable $\ka$. 
\csju

\noi{\bf Note 2.} If $M$ is not $\ka$-slender, then, every morphism $f:\prod_I A_i\lra M$, such that for $J=\sst{i\in I}{f|A_i\neq 0}$, $|J|\geq\ka$ will be called a {\it non-slender morphism}. Given a non-$\ka$-slender module $M$, such a non-slender morphism always exists and, by taking the appropriate restriction to $\prod_JA_i$ we may then assume that, for a non-$\ka$-slender module $M$, there is a morphism $f:\prod_IA_i\lra M$ with $|I|=\ka$, such that $f|A_i\neq 0$, for every $i\in I$. 
\csju

We note immediately that, in the definition, we may take any index set of cardinality $>\ka$ as well as that we may replace all $A_i$ by cyclic modules or by  the identical objects, namely the ground ring $R$ as detailed in the following: 

\proclaim{Theorem 3} Given an infinite cardinal $\ka$ and an $M\in\rmod$, the following are equivalent:
\roster
\item $M$ is $\ka$-slender.
\item $\forall I$\,\, $|I|\geq\ka$\,\, $\forall A_i\in\rmod, i\in I$, for every morphism $f:\prod_{i\in I}A_i\lra M$, $|\sst{i\in I}{f|A_i\neq 0}|<\ka.$
\item $\forall I$\,\, $|I|\geq\ka$, for every morphism $f:\prod_{i\in I}R_i\lra M, \forall i R_i=R$,
  \newline  $|\sst{i\in I}{f(\e_i)\neq 0}|<\ka.$
\item $\forall I$\,\, $|I|=\ka$, for every morphism $f:\prod_{i\in I}R_i\lra M, \forall i R_i=R$,
  \newline  $|\sst{i\in I}{f(\e_i)\neq 0}|<\ka.$
\item $\forall I$\,\, $|I|\geq\ka$, for every morphism $f:\prod_{i\in I}Ra_i\lra M$, 
  \newline  $|\sst{i\in I}{f(a_i)\neq 0}|<\ka.$
\item $\forall I$\,\, $|I|=\ka$, for every morphism $f:\prod_{i\in I}Ra_i\lra M$, 
  \newline  $|\sst{i\in I}{f(a_i)\neq 0}|<\ka.$
\endroster
In case of regular $\ka$, we also have the following equivalent statements:
\roster
\item "(7)" $\forall I$\,\, $|I|=\ka$\,\, $\forall A_i\in\rmod, i\in I$, for every morphism $f:\prod_{i\in I}A_i\lra M$, $\exists i_0<\ka$ such that  $\forall i>i_0\,\, f|A_i=0$.
\item "(8)" $\forall I$\,\, $|I|=\ka$,\,  for every morphism $f:\prod_{i\in I}Ra_i\lra M$, $\exists i_0<\ka$ such that  $\forall i>i_0\,\, f(a_i)=0$.
\item "(9)" $\forall I$\,\, $|I|=\ka$,\,  for every morphism $f:\prod_{i\in I}R_i\lra M$, $\forall i R_i=R$,
$\exists i_0<\ka$ such that  $\forall i>i_0\,\, f(\e_i)=0$.
\item "(10)" $\forall I$\,\, $cf\,|I|\geq\ka$\,\, $\forall A_i\in\rmod, i\in I$, for every morphism $f:\prod_{i\in I}A_i\lra M$, $\exists i_0<\ka$ such that  $\forall i>i_0\,\, f|A_i=0$.
\item "(11)" $\forall I$\,\, $cf\,|I|\geq\ka$,\,  for every morphism $f:\prod_{i\in I}Ra_i\lra M$, $\exists i_0<\ka$ such that  $\forall i>i_0\,\, f(a_i)=0$.
\item "(12)" $\forall I$\,\, $cf\,|I|\geq\ka$,\,  for every morphism $f:\prod_{i\in I}R_i\lra M$, $\forall i R_i=R$,
$\exists i_0<\ka$ such that  $\forall i>i_0\,\, f(\e_i)=0$.
\endroster
\endproclaim 
\prf (1)$\shimply (2)$: Let $|I|>\ka$ and $f:\prod_{i\in I}A_i\lra M$. If, on the contrary, $\exists J\se I$, $|J|=\ka$ such that $\forall j\in J$\,\, $f|A_j\neq 0$, then we can take the restriction $f\pr=f|\prod_{j\in J}A_j$ with the same coordinate property. This would then contradict (1). 
(2)$\shimply$ (3),  (3)$\shimply$ (4), (5) $\shimply$ (6) hold because, respectively,   (2) is nominally more general than (3) and (3) is nominally more general than (4) just as (5) is nominally more general than (6). 
(4)$\shimply$(5): Let, on the contrary, $\exists J\se I$, $|J|=\ka$, such that $\forall j\in J$, $f(a_j)\neq 0$. We have the quotient maps $q_j:R\lra R/\ann(a_j)\cong Ra_j$ and the product map $q=\prod q_j:\prod_{j\in J} R_j\lra \prod_{j\in J}Ra_j$. Consider $f\pr=fq:\prod_{j\in J}R_j\lra M$. We have $f\pr(\e_j)=f(a_j)\neq 0$, $\forall j\in J$, which would contradict (4). 
(6)$\shimply$ (1): Let $|I|=\ka$ and $f:\prod_{i\in I}A_i\lra M$ be such that, on the contrary,  $\exists J\se I$\,\, $|J|=\ka$ with $\forall j\in J\,\, f|A_j\neq 0$; in other words $\exists a_j\in A_j$ with $f|Ra_j\neq 0$ 
Consider the restriction of $f$, namely $f\pr:\prod_{j\in J}Ra_j\lra M$. But then $f\pr|Ra_j$ on a set $J$ of cardinality $\ka$ which would contradict (6). 
The equivalences (7) -- (12) are proved in a like manner as (1) -- (6). We only need to connect the two batches: (7)$\shimply$(1) follows, once we note that $\forall i_0<\ka$\,\, $|(\lla, i_0)|<\ka$. As for (1)$\shimply$(7), given a morphism  $f:\prod_{i\in I}A_i\lra M$, the cardinality $|S_i|=|\sst{i<\ka}{f|A_i\neq 0}|<\ka$. This implies that sup$S_i=i_0<\ka$, since $\ka$ was assumed to be regular.
We note that regularity of $\ka$ is not needed for implications from the second batch of statements to the first. 
\qed

It appears that $\ka$-slenderness is a characteristic of the lattice of submodules as may be seen from the following:

\proclaim{Proposition 4} 
\roster
\item The trivial module $0$ is $\ka$-slender, for every $\ka$. 
\item $M$ is $\ka$-slender, iff $\forall N\leq M$,\,\, $N$ is $\ka$-slender. 
\item For every infinite cardianl $\ka$, every slender $R$-module is $\ka$-slender.
\item Let $\ka\leq\lambda$; then every $\ka$-slender module is also $\lambda$-slender. 
\item $R^\ka$ is not $\ka$-slender.
\item For all cardinals $\ka<\lambda$, and $B_j\in\rmod, j<\ka$, $\prod_{j<\ka}B_j$ is $\lambda$-slender, if and only if every $B_j$ is $\lambda$-slender. In particular, 
$R^\ka$ is $\lambda$-slender if and only if   $R$ is $\lambda$-slender. Furthermore, $\oplus_{j<\ka}B_j$ is $\lambda$-slender if and only if every $B_j$ is $\lambda$-slender. 
\endroster
\endproclaim

\prf (1) The definition verifies trivially. 

(2)  If $f:R^\ka\lra N\hookr M$, then use Theorem 3(4) to conclude that $\sst{i\in I}{f(\e_i)\neq 0}|<\ka$, which establishes $\ka$-slenderness of $N$. The other direction is a tautology.

(3) By a known result (see e.g. Dimitric (2017), Theorem 3.10), an object $M$ is slender iff for every index set $I$ and every morphism $f:\prod_{i\in I}A_i\lra M$,  $\sst{i\in I}{f|A_i\neq 0}$ is finite (hence $<\ka$).

(4) Let $M$ be $\ka$-slender and let $f:\prod_{i\in I}Ra_i\lra M$, ($|I|=\lambda$). By Theorem 3(5) $|\sst{i\in I}{f(a_i)\neq 0}|<\ka<\lambda$, which establishes $\lambda$-slenderness of $M$.

(5) The identity map $id:R^\ka\lra \prod_{j<\ka}R_j$ is such that $|\sst{i\in I}{id(\e_i)\neq 0}|=\ka$, which shows, by Theorem 3, that $R^\ka$ is not $\ka$-slender. 

(6)  If $\prod_{j<\ka}B_j$ is $\lambda$-slender, then, by (2), every submodule is $\lambda$-slender, hence that applies to each $B_j$ as well. Now assume that every $B_j$ is $\lambda$-slender. If $f:\prod_{i<\lambda}Ra_i\lra \prod_{j<\ka}B_j$, we know that then $f=\left(f_j:\prod_{i<\lambda}Ra_i\lra B_j\right)_{j\in J}$, $|J|=\ka$. We know that $\forall j\in J$, for\,$S_j=\sst{i\in I}{f_j(a_i)\neq 0}$, $|S_j|<\lambda$, since every $B_j$ is $\lambda$-slender. We have $\sst{i\in I}{f(a_i)\neq 0}=\cup_{j\in J}S_j$. Assume first that $\lambda$ is a regular cardinal. Then $|\cup_{j\in J}S_j|\leq\sum_{j\in J}|S_j|<\ka\lambda=\lambda$ (the latter strict inequality holds because $\lambda$ is regular). Thus, indeed $\prod_{j<\ka}B_j$ is $\lambda$-slender.  Hence, this statement is true for regular cardinal $\lambda=\ka^+$, namely $\prod_{j<\ka}B_j$ is $\ka^+$-slender. By (4), $\prod_{j<\ka}B_j$ is $\lambda$-slender, for every $\lambda\geq\ka^+>\ka$. 
The remaining claims are a special case and the fact that the direct sum is submodule of the direct product, hence by (2) has to be slender. \qed

\csju

 Given an infinite cardinal $\ka$, then a submodule $N\leq M\in\rmod$ is said to be $\ka$-{\it pure}
in $M$, if every system of ($|I|<\ka$) equations of the form 
$$\sum_{j\in J}r_{ij}x_j=n_i\in N,\quad i\in I,
\quad r_{ij}\in R, \eqno{(1)}
$$ 
with $<\ka$ unknowns $x_j, j\in J, |J|<\ka$
that has a solution in $M^J$, also has a solution in $N^J$. Notation for this is
$N\leq_{\ka*} M$. Thus purity is then same as $\aleph_0$-purity. The derivative notion of a $\ka$-pure exact sequence is straightforward.  
As $\ka$ is increased, the classes (sets) of $\ka$-pure exact sequences
get smaller, in general. 
A module is $\ka$-pure injective if it has injective property with respect to all $\ka$-pure exact sequences. A module $M$ is equationally (algebraically) $\ka$-{\it compact}, if every system of $\leq\ka$ linear equations:
$$\sum_{j\in J}r_{ij}x_j=m_i\in M,\quad i\in I,\quad r_{ij}\in
R\quad \eqno (2)
$$ 
with the property that every finite subsystem has a solution, then has a global solution. A module is algebraically compact iff it is $\ka$-compact, for every
cardinal $\ka$. Given $\ka<|R|$ one can construct examples of $\ka$-compact
modules that are not algebraically compact. However, if 
 $M\in\rmod$ is $\ka$-algebraically compact, for some $\ka\geq
|R|$, then $M$ is algebraically compact.

We have mimicked \L o\'s (1959) to produce the following result, needed in the sequel:

\proclaim{Theorem 5} Let $A$ be an index set,  
$\Cal F$  a $\ka$-complete filter on $A$ and $M_\al\in\rmod$, $\al\in A$; 
then $\prod^{\Cal F}_{\al\in A} M_\al$ is $\ka$-pure in 
$\prod_{\al\in A} M_\al$.
Specially, for $\Cal F=\Cal F_0$, the coproduct $\op M_\al$ is pure in $\prod M_\al$.
\endproclaim
\prf Assume that the system of linear equations 
$$\sum_{j\in J}r_{ij}x_j=n_i\in \mathop{{\prod_{\al\in A}}^{\Cal F}}M_\al,\quad i\in I,
\quad (r_{ij})_{I\times J} \text{ row finite },\quad |I|,|J|<\ka
$$ 
has a solution $m_j\in\prod_{\al\in A}M_\al, j\in J$; this then translates into
the componentwise equalities:
$\sum_{j\in J}r_{ij}m_{j\al}=n_{i\al}\,\,\, (*)$, $\al\in A$. By definition $zero(n_i)\in\Cal F$,
for all $i\in I$, and since $|I|<\ka$, we get, by $\ka$-completeness of $\Cal F$, 
that $Z=\cap_{i\in I}zero(n_i)\in\Cal F$. Now define $y_j\in\prod_{\al\in I}M_\al$
componentwise: $y_{j\al}=m_{j\al}$, if $\al\notin Z$ and $y_{j\al}=0$, if $\al\in Z$. 
Every $zero(y_j)\spe Z\in\Cal F$, thus all $y_j\in\prod^{\Cal F} M_\al$; but the $y_j$
also provide a solution of the original system of equations, by the way we defined
them, by ($*$) and by the fact that for $\al\in Z$ we have $n_{i\al}=0$, for all $i$.
\qed
\csju

Denote by $\Cal S_\ka$ the class of $\ka$-slender modules, where $\Cal S$ denotes, for brevity, the class of slender modules. 

\proclaim{Proposition 6}
\roster
\item We have an ascending chain of non-empty classes: 
$$
\Cal S\subset\dots\subset\Cal S_\ka\subset \Cal S_{\ka^+}\subset\dots\subset\Cal S_\lambda\subset\dots\subset\rmod\quad \ka<\lambda\,.
$$
The chain is strictly ascending, if $R$ is slender. 
\item The union of this chain is $\neq\rmod$, since non-zero algebraically compact modules are not $\ka$-slender, for any $\ka$.
\endroster
\endproclaim
\prf (1) is a consequence of Proposition 4(1),(4). For (2), given a cardinal $\ka$, assume that  $M\in \rmod$ is algebraically compact and let $0\neq a\in M$. By Theorem 5, we have a pure exact sequence 
$$ 0\lra \bigoplus_I R_i\lra\prod_IR_i\lra\prod_IR_i/\bigoplus_IR_i\lra 0, \quad R_i=R  \eqno (3)
$$
Define $f_0:\oplus R_i\lra M$ coordinatewise: $\forall i\in I f_0(\e_i)=a$.
Since $M$ is algebraically compact, we can extend $f_0$ to the morphism $f:\prod_IR_i\lra M$, for which we have $\forall i, f(\e_i)=a\neq 0$, which shows that $M$ is not $\ka$-slender, for any $\ka$. \qed

Consequently, if $R$ algebraically compact (pure injective), then, by Proposition 6(2), $R$  is not $\lambda$-slender, for any $\lambda$ and then the product $R^\ka$, being algebraically compact, is not $\lambda$-slender, for any $\lambda, \ka$. 

\vfill
\eject

\noi{\bf 2. Tailwise slenderness}
\csju
\noi{\bf Definition 7.} Given a cardinal $\ka$, an $M\in\rmod$ is said to be {\it tailwise $\ka$-slender}, or {\it t$\ka$-slender} for short, if for every morphism $f:\prod_{i<\ka}Ra_i\lra M$, there exists an $i_0<\ka$ such that 
$f(\prod_{i\geq i_0}Ra_i)=0$. This is equivalent to the requirement that, for every morphism $f:\prod_{i<\ka}R_i\lra M$, $R_i=R$, there exists an $i_0<\ka$ such that $f(\prod_{i\geq i_0}R_i)=0$.

We note a straightforward but important fact as follows:

\proclaim{Proposition 8}
\roster
\item If $M$ is $t\ka$-slender, then it is $\ka$-slender. 
\item If $M$ is t$\ka$-slender, then, for all cyclic modules $Ra_i, i\in I$, $|I|=\ka$: 
$$\homc{R}{\prod_{i\in I}Ra_i/\mathop{{\prod_{i\in I}}^{\ka}}Ra_i}{M}=0.$$
\endroster
\endproclaim
\prf (1) $S_i=\sst{i<\ka}{f|A_i\neq 0}\se(\lla, i_0)$ and $|(\lla, i_0)|<\ka$, since $\ka$ is a cardinal. 

(2) If $\ka=|I|$ and $f:\prod_{i\in I} Ra_i\lra M$ is a morphism, then, by   t$\ka$-slenderness of  $M$, $\exists i_0<\ka$ such that $f(\prod_{i\geq i_0}Ra_i)=0$. We note that $\prod_{i<i_0}Ra_i\se\prod_{i\in I}^\ka Ra_i$ since $\ka$ is a cardinal and $|(\lla,i_0)|<\ka$. The claim will follow, once we note the obvious splitting: $\prod_{i<\ka}Ra_i=\prod_{i<i_0}Ra_i
\oplus \prod_{i\geq i_0}Ra_i$.
\qed
 \csju

\noi As for the converse of implication (1) in  this proposition, it may not always be true and it is related to intricate constructions of set-theoretical nature. It is well-known that, for $\ka=\om_0$, the equivalence holds, if and only $\ka$ is a non-measurable cardinal (see Dimitric (2017), Theorem 3.10).

\proclaim{Proposition 9} 
\roster
\item The trivial module $0$ is t$\ka$-slender, for every $\ka$. 
\item $M$ is t$\ka$-slender, iff $\forall N\leq M$,\,\, $N$ is t$\ka$-slender. 
\item For every infinite non-measurable cardianl $\ka$, every slender $R$-module is t$\ka$-slender.
\item $R^\ka$ is not t$\ka$-slender.
\item For all cardinals $\ka< cf\lambda$, and $\sst{B_j\in\rmod}{j<\ka}$, $\prod_{j<\ka}B_j$ is t$\lambda$-slender, if and only if every $B_j$ is t$\lambda$-slender. In particular, 
$R^\ka$ is t$\lambda$-slender if and only if   $R$ is t$\lambda$-slender. Furthermore, $\oplus_{j<\ka}B_j$ is t$\lambda$-slender if and only if every $B_j$ is t$\lambda$-slender. 
\item If $0\lra A\ese{\al}B\ese{\be}C\lra 0$ is an exact sequence and $A, C$ are t$\ka$-slender, then $B$ is likewise t$\ka$-slender. 
\endroster
\endproclaim

\prf (1) and (2) follow directly from the definition. 
(3) follows from Dimitric (2017), Theorem 3.10(5), which  states, that if $|I|=\ka$ is a non-measurable cardinal and $M$ is slender, then, for every morphism $f:\prod_{i\in I}A_i\lra M$, there exists an $i_0\in I$, $i_0<\om$ such that $f(\prod_{i\geq i_0}A_i)=0$. For (4), consider the non-t$\ka$-slender identity morphism $id: R^\ka\lra R^\ka$. 

(5) If $\prod_{j<\ka}B_j$ is t$\lambda$-slender, then, by (2), every submodule is t$\lambda$-slender, hence that applies to each $B_j$ as well. Now assume that every $B_j$ is t$\lambda$-slender. If $f:\prod_{i<\lambda}Ra_i\lra \prod_{j<\ka}B_j$, we know that then $f=\left(f_j:\prod_{i<\lambda}Ra_i\lra B_j\right)_{j\in J}$, $|J|=\ka$. We have $\forall j\in J$ $\exists i(j)<\lambda$ with $f(\prod_{i\geq i(j)}Ra_i)=0$. Let $i_0=\text{sup }\sst{i(j)}{j<\ka}$. By the assumption, cf$\lambda>\ka$, therefore $i_0<\lambda$. Now we clearly have  $f(\prod_{i\geq i_0}Ra_i)=0$. 
The remaining claims are a special case and the fact that the direct sum is submodule of the direct product, hence by (2) has to be t$\lambda$-slender. 

(6) Let $f:\prod_{i<\ka}Ra_i\lra B$ be an arbitrary morphism. Then $\be f:\prod_{i<\ka}Ra_i\lra C$, hence, by t$\ka$-slenderness of $C$, there is an $i\pr<\ka$ such that $\be f(\prod_{\ka>i\geq i\pr}Ra_i)=0$. In other words,
$f(\prod_{\ka>i\geq i\pr}Ra_i)\se\ker\be=\ima\al\cong A$, which implies that $f$ maps $\prod_{\ka>i\geq i\pr}Ra_i$ into $\ima\al\cong A$. On the other hand, $A$ was assumed to be t$\ka$-slender which then implies that there is an $i_0<\ka$, $i_0\geq i\pr$ with $f(\prod_{\ka>i\geq i_0}Ra_i)=0$, which establishes slenderness of $B$.
\qed

\ju
\noi{\bf 3. Classes $\Cal H_\ka$}
\csju

\noi{\bf Definition 10.} Given an infinite cardinal $\ka$, an $M\in\rmod$ is called an $h\ka$-{\it module}, if, for every index set $I$, and every family of $R$-modules $\sst{A_i}{i\in I}$, the following holds:
$$\homc{R}{\prod_{i\in I}A_i/\mathop{{\prod_{i\in I}}^{\ka}}A_i}{M}=0.$$
For brevity, denote 
$D_\ka=\prod_{i\in I}A_i/\prod_{i\in I}^{\,\,\ka}A_i$, so that we can rewrite this condition as 
$\homc{R}{D_\ka}{M}=0.$

The class of $h\ka$-modules is denoted by $\Cal H_\ka$.

Some well-known properties of the Hom functor are instrumental in obtaining some properties of $h\ka$ modules as follows:

\proclaim{Proposition 11}
\roster
\item $0\in\Cal H_\ka$.
\item $\Cal H_k$ is closed with respect to submodules. 
\item $\prod M_j\in \Cal H_k$ if and only if, every $M_j\in\Cal H_\ka$ (closure with respect to products). 

\item For a short exact sequence $0\lra A\ese{\al}B\ese{\be}C\lra 0$,  if $A, C\in \Cal H_\ka$, then $B\in\Cal H_\ka$. 
\item If $M$ is $t\ka$-slender, then $M\in\Cal H_\ka$.
\item If $\ka<\lambda$, then $\Cal H_\ka\se\Cal H_{\lambda}$.
\endroster
\endproclaim

\prf 
 (1) is trivial. (2) follows from the fact that, if $B\leq A$, then $\hom{D_\ka}{B}\leq\hom{D_\ka}{A}$. The natural isomorphism $\hom{D_\ka}{\prod M_j}\cong\prod\hom{D_\ka}{M_j}$ establishes one implication of (3) and the other implication is a consequence of (2). (4)	 is a consequence of left exactness of the Hom functor: $0\lra \hom{D_\ka}{A}\ese{}\hom{D_\ka}{B}\ese{}\hom{D_\ka}{C}$.  For (5) assume that, on the contrary,  there were an $a=(a_i)_{i\in I}\in\prod_{i\in I}A_i$, and a morphism $f:D_\ka\lra M$ with $f(\ov a)=0$; then it would contradict Proposition 8(2) since we would have 
$\homc{R}{\prod_{i\in I}Ra_i/{\prod_{i\in I}}^{\ka}Ra_i}{M}\neq 0.$
For (6), use the fact that, for $\ka<\lambda$, 
${\prod_{i\in I}}^{\ka}A_i\leq{\prod_{i\in I}}^{\lambda}A_i$.
\qed
\sju

Properties (2)--(4) signify that $\Cal H_\ka$ is a torsion free class for a torsion theory, for every $\ka$ (cf. e.g. 
Stenstr\"om (1975)). 
\sju

A good question is whether Proposition 11(5) holds for coordinate $\ka$-slender modules. It does for $\ka=\om_0$ and non-measurable index sets $I$ (cf. Dimitric (2017), Theorem 3.9). We are exploring this issue for uncountable $\ka$, at present time. 

\ju
{\bf References}
\ju

Dimitric, R. (2017):
{\it Slender Modules and Rings. Vol. I. Foundations}. Cambridge: Cambridge University Press. To appear. 

\L o\'s, Jerzy (1959):
Linear equations and pure subgroups, {Bulletin de l'Acad\'emie 
Polonaise des Sci., S\'erie math.}, {\bf 7}, No.1, 13--18.

Stenstr\"om, Bo (1975): {\it Rings of Quotients. An Introduction to Methods of Ring Theory}. Berlin, Heidelberg, New York: Springer-Verlag.

\ju

May 25, 2017

\copyright \,\, R. Dimitric  1985, 2017

\end
\enddocument